\documentclass[11pt]{article}
\usepackage[T1]{fontenc}
\usepackage{amsmath,amssymb,amsthm}
\usepackage{graphicx}
\usepackage{subcaption}
\usepackage{tikz}
\usetikzlibrary{shapes.misc, shapes.geometric, positioning, automata}
\usepackage{booktabs}
\usepackage[margin=1in]{geometry}
\usepackage{xcolor}
\usepackage{enumitem}
\usepackage{fancyhdr}
\usepackage{mathtools}
\usepackage{url}
\usepackage{float}
\usepackage[colorlinks=true, linkcolor=blue!70!black, citecolor=green!60!black, urlcolor=red!70!black]{hyperref}

\definecolor{theoremblue}{RGB}{0, 102, 204}
\definecolor{lemmagreen}{RGB}{0, 153, 51}
\definecolor{proofgray}{RGB}{245, 245, 245}
\definecolor{remarkred}{RGB}{204, 0, 0}
\definecolor{definitionpurple}{RGB}{102, 0, 153}

\usepackage[theorems,skins,breakable]{tcolorbox}

\newtcbtheorem[number within=section]{mythm}{Theorem}{
  enhanced,
  breakable,
  colback=blue!5,
  colframe=theoremblue,
  coltitle=white,
  fonttitle=\bfseries,
  sharp corners,
  attach boxed title to top left={yshift=-2mm,xshift=5mm},
  boxed title style={colback=theoremblue,sharp corners}
}{thm}

\newtcbtheorem[use counter from=mythm]{mylem}{Lemma}{
  enhanced,
  breakable,
  colback=green!5,
  colframe=lemmagreen,
  coltitle=white,
  fonttitle=\bfseries,
  sharp corners,
  attach boxed title to top left={yshift=-2mm,xshift=5mm},
  boxed title style={colback=lemmagreen,sharp corners}
}{lem}

\newtcbtheorem[use counter from=mythm]{myconj}{Conjecture}{
  enhanced,
  breakable,
  colback=yellow!10,
  colframe=yellow!50!black,
  coltitle=white,
  fonttitle=\bfseries,
  sharp corners,
  attach boxed title to top left={yshift=-2mm,xshift=5mm},
  boxed title style={colback=yellow!60!black,sharp corners}
}{conj}

\newtcbtheorem[number within=section]{mydef}{Definition}{
  enhanced,
  breakable,
  colback=purple!5,
  colframe=definitionpurple,
  coltitle=white,
  fonttitle=\bfseries,
  sharp corners,
  attach boxed title to top left={yshift=-2mm,xshift=5mm},
  boxed title style={colback=definitionpurple,sharp corners}
}{def}

\newtcbtheorem[use counter from=mythm]{myprop}{Proposition}{
  enhanced,
  breakable,
  colback=orange!5,
  colframe=orange!80!black,
  coltitle=white,
  fonttitle=\bfseries,
  sharp corners,
  attach boxed title to top left={yshift=-2mm,xshift=5mm},
  boxed title style={colback=orange!80!black,sharp corners}
}{prop}

\newtcbtheorem[use counter from=mythm]{mycor}{Corollary}{
  enhanced,
  breakable,
  colback=cyan!5,
  colframe=cyan!70!black,
  coltitle=white,
  fonttitle=\bfseries,
  sharp corners,
  attach boxed title to top left={yshift=-2mm,xshift=5mm},
  boxed title style={colback=cyan!70!black,sharp corners}
}{cor}

\newtcbtheorem[use counter from=mythm]{mythmruskey}{Theorem (Ruskey \& Deugau)}{
  enhanced,
  breakable,
  colback=blue!5,
  colframe=theoremblue,
  coltitle=white,
  fonttitle=\bfseries,
  sharp corners,
  attach boxed title to top left={yshift=-2mm,xshift=5mm},
  boxed title style={colback=theoremblue,sharp corners}
}{ruskey}

\newtcolorbox{remark}[1][]{
  enhanced,
  breakable,
  colback=red!5,
  colframe=remarkred,
  coltitle=remarkred,
  title={\textbf{Remark}#1},
  sharp corners,
  boxrule=1pt,
  left=5pt,right=5pt,top=5pt,bottom=5pt
}

\renewenvironment{proof}[1][\proofname]{%
  \begin{tcolorbox}[
    enhanced,
    breakable,
    colback=proofgray,
    colframe=gray!50,
    coltitle=black,
    title={\textbf{#1}},
    sharp corners,
    boxrule=0.5pt,
    left=5pt,right=5pt,top=5pt,bottom=5pt
  ]
}{%
  \hfill$\blacksquare$
  \end{tcolorbox}
}

\newcommand{\N}{\mathbb{N}}
\newcommand{\Z}{\mathbb{Z}}
\newcommand{\floor}[1]{\lfloor #1 \rfloor}

\newcommand{\oeis}[1]{\href{https://oeis.org/#1}{\texttt{#1}}}

\title{\Large\textbf{A study of a family of self-referential sequences:} \\[0.5em] 
\large Asymptotics, exact forms and combinatorial connections}
\author{Benoit Cloitre}
\date{\today}

\begin{document}

\maketitle

\begin{abstract}
We introduce and analyze a three-parameter family of self-referential integer sequences $S(x,y,z)$: starting from $a(1)=x$, each term advances by $y$ when the index $k$ has already appeared as a value and by $z$ otherwise. This simple rule generates a surprising zoo of behaviors, many of which are catalogued—albeit in a rather unstructured fashion—in the OEIS. This family has recently and independently been studied by Fokkink and Joshi, who named them "hiccup sequences" and established their general morphic nature. Our work provides a complementary, in-depth analysis of major subfamilies. Whenever $y>z>0$, we prove that the density $a(k)/k$ converges to the positive root of $r^{2}-zr-(y-z)=0$. Two subfamilies, $S(x,Z+1,Z)$ and $S(x,Z,Z+1)$, yield explicit non-homogeneous Beatty sequences, providing explicit formulas for numerous OEIS entries. For $y=0$ and $z\ge 2$, the sequences eventually become periodic and satisfy linear recurrences. Critical cases with a zero discriminant unveil geometric patterns on triangular, square, and hexagonal lattices. Finally, via tree-like representations we uncover a tight link with meta-Fibonacci recurrences. These results position $S(x,y,z)$ as a unifying framework connecting additive combinatorics, number theory, and discrete dynamics.
\end{abstract}

\section{Introduction and definitions}\label{sec:intro}

Self-referential integer sequences, where the rule for generating the next term depends on the sequence's history, often yield intricate global structures from simple local rules. Their appeal lies in understanding the emergence of complex, predictable behavior from recursive definitions. Such sequences frequently appear in combinatorial number theory and have connections to fractal geometry, dynamical systems, and combinatorics. In our 2003 paper, "Numerical analogues of Aronson's sequence" \cite{CloitreSloaneVandermast2003}, we announced our intention to study the following family of sequences.
\begin{mydef}{The $S(x,y,z)$ family}{sxyz}
Let $x,y,z\in\N$. The sequence $(a(k))_{k\ge1}$ is defined by:
\begin{enumerate}[label=(\roman*)]
\item $a(1)=x$.
\item For $k>1$, let $\mathcal{A}_{k-1}=\{a(1),\dots,a(k-1)\}$. Then,
\[
a(k)=\begin{cases}
a(k-1)+y & \text{if }k\in\mathcal{A}_{k-1}\quad(\text{a \emph{hit}}),\\
a(k-1)+z & \text{if }k\notin\mathcal{A}_{k-1}\quad(\text{a \emph{miss}}).
\end{cases}
\]
\end{enumerate}
\end{mydef}

\begin{remark}[Recent related work]
This family of sequences was recently and independently investigated by Fokkink and Joshi \cite{fokkink cloitre xyz.tex}. They introduce the slightly more general notion of a $(j,x,y,z)$-"hiccup sequence", where the condition depends on $n-j \in \mathcal{A}_{n-1}$. Our family $S(x,y,z)$ corresponds precisely to their $(0,x,y,z)$-hiccup sequences. A major result of their work, which extends previous results from Bosma, Dekking, and Steiner, is a general proof that the characteristic sequences of all such "hiccup sequences" are morphic.
\end{remark}

At the time of our initial work, a lack of a cohesive framework and general results compelled us to postpone a detailed study. Over twenty years later, this paper fulfills that original goal. While the work of Fokkink and Joshi provides a broad framework establishing the morphic nature of these sequences, our focus here is a deep-dive into the rich analytic and combinatorial properties of key sub-families, most of which are not explored elsewhere. We believe this work demonstrates the rich structure of these sequences, offering significant points of interest for number theorists and combinatorists alike, particularly in the areas of combinatorics on words and enumerative problems.
The structure of this paper is as follows. In Section~\ref{sec:asymptotic}, we establish the general asymptotic linearity of the sequences for the case $y > z > 0$. Section~\ref{sec:beatty} is dedicated to two important subfamilies that generate non-homogeneous Beatty sequences, providing explicit formulas for many OEIS entries. In Section~\ref{sec:recurrence}, we analyze the case where $y=0$, which leads to ultimately periodic sequences satisfying linear recurrences. Section~\ref{sec:vanishing} explores cases where the characteristic discriminant vanishes, revealing deep connections to combinatorial objects like triangular, square, and hexagonal numbers. Finally, in Section~\ref{sec:metafib}, we establish a formal bridge between our sequences and the theory of meta-Fibonacci sequences, uncovering a shared combinatorial foundation.

\section{Asymptotic linearity for \texorpdfstring{$y>z>0$}{y>z>0}}\label{sec:asymptotic}

We first establish that when the hit increment $y$ exceeds the miss increment $z$, the sequence exhibits a stable asymptotic slope. Let $r_{0}$ be the positive root of the characteristic equation $r^{2}-zr-(y-z)=0$:
\[
r_{0}=\frac{z+\sqrt{z^{2}+4(y-z)}}{2}.
\]
This root satisfies the equilibrium equation $r_{0}=z+\frac{y-z}{r_{0}}$.
\begin{mythm}{Main asymptotic result}{main}
For a sequence $S(x,y,z)$ with $y>z>0$,
\[
\lim_{n\to\infty}\frac{a(n)}{n}=r_{0}.
\]
\end{mythm}

\begin{proof}
\textbf{Step 1: Exact relation.}
By definition:
\[
a(n) - a(n-1) = z + (y - z) \cdot \mathbf{1}_{\{n \in \mathcal{A}_{n-1}\}}.
\]
Summing for $k = 2$ to $n$:
\[
a(n) - a(1) = z(n-1) + (y - z) H_n,
\]
where $H_n = \sum_{k=2}^n \mathbf{1}_{\{k \in \mathcal{A}_{k-1}\}}$ is the number of hits up to step $n$. Thus:
\[
\boxed{a(n) = x + z(n-1) + (y - z) H_n.}
\]

\textbf{Step 2: Link between $H_n$ and the counting function.}
Let $N(t) = \max\{k \mid a(k) \leq t\}$ be the counting function of the sequence.
\begin{itemize}[leftmargin=*]
\item \textbf{Key property:} $k \in \mathcal{A}_{k-1}$ if and only if $k$ is in the image of the sequence, i.e., $\exists j < k$ such that $a(j) = k$.
\item Thus $H_n = |\{k \in \{2,\dots,n\} \mid k \in \operatorname{im}(a)\}|$.
\item Since $N(n) = \max\{k \mid a(k) \leq n\}$, we have:
\[
H_n = N(n) - \mathbf{1}_{\{1 \in \operatorname{im}(a)\}} + O(1) = N(n) + O(1),
\]
\end{itemize}
since $\mathbf{1}_{\{1 \in \operatorname{im}(a)\}}$ is constant (0 or 1). We deduce:
\[
\frac{H_n}{n} = \frac{N(n)}{n} + O\left(\frac{1}{n}\right).
\]

\textbf{Step 3: Asymptotic analysis via $\limsup$ and $\liminf$.}
Let:
\[
\overline{r} = \limsup_{n \to \infty} \frac{a(n)}{n}, \quad \underline{r} = \liminf_{n \to \infty} \frac{a(n)}{n}.
\]
From the equation in Step 1:
\[
\frac{a(n)}{n} = \frac{x}{n} + z \frac{n-1}{n} + (y - z) \frac{H_n}{n}.
\]
Taking upper and lower limits, and using $\frac{H_n}{n} = \frac{N(n)}{n} + O(1/n)$:
\[
\boxed{\overline{r} = z + (y - z) \limsup_{n \to \infty} \frac{N(n)}{n}, \quad \underline{r} = z + (y - z) \liminf_{n \to \infty} \frac{N(n)}{n}.}
\]

\textbf{Step 4: Duality between sequence and counting function.}
\begin{mylem}{Density-counting duality relation}{density}
Let $(a(n))_{n \in \mathbb{N}^*}$ be a strictly increasing sequence of real numbers such that $\lim_{n \to \infty} a(n) = +\infty$. Let $N: \mathbb{R} \to \mathbb{N}$ be the counting function defined by $N(t) = \max\{k \in \mathbb{N} \mid a(k) \leq t\}$. Define the limit superior and limit inferior of the sequence's rate of growth as:
\[
\overline{r} = \limsup_{n \to \infty} \frac{a(n)}{n} \quad \text{and} \quad \underline{r} = \liminf_{n \to \infty} \frac{a(n)}{n}.
\]
Then the following duality relations hold:
\[
\limsup_{t \to \infty} \frac{N(t)}{t} = \frac{1}{\underline{r}} \quad \text{and} \quad \liminf_{t \to \infty} \frac{N(t)}{t} = \frac{1}{\overline{r}}.
\]
\end{mylem}

\paragraph{Proof of Lemma~\ref{lem:density}.}
The proof is based on the fundamental inequality derived from the definition of the counting function $N(t)$.
\textbf{Step 1: The fundamental inequality.}
By definition, $a(N(t))$ is the largest term in the sequence that is less than or equal to $t$. This directly implies:
\[
a(N(t)) \leq t < a(N(t)+1).
\]
This inequality holds for all $t \geq a(1)$.
\textbf{Step 2: Bounding the ratio $N(t)/t$.}
We manipulate the fundamental inequality to establish a tight bound on the ratio $N(t)/t$. From the left side, $a(N(t)) \leq t$, we invert to get $\frac{1}{t} \leq \frac{1}{a(N(t))}$, which implies:
\[
\frac{N(t)}{t} \leq \frac{N(t)}{a(N(t))}.
\]
From the right side, $t < a(N(t)+1)$, we similarly get $\frac{1}{a(N(t)+1)} < \frac{1}{t}$, which implies:
\[
\frac{N(t)}{a(N(t)+1)} < \frac{N(t)}{t}.
\]
Combining these gives the sandwich inequality:
\[
\frac{N(t)}{a(N(t)+1)} < \frac{N(t)}{t} \leq \frac{N(t)}{a(N(t))}.
\quad (*)
\]
\textbf{Step 3: Asymptotic analysis.}
We now analyze the limits of the bounding terms in $(*)$ as $t \to \infty$. Let $v_k = k/a(k)$. We are comparing the asymptotic behavior of $v_{N(t)}$ with that of the sequence $(v_k)$. The function $N(t)$ is a non-decreasing step function from $[a(1), \infty)$ onto $\mathbb{N}$. Because it is non-decreasing and surjective onto the integers, the set of accumulation points of the function $v_{N(t)}$ as $t \to \infty$ is identical to the set of accumulation points of the sequence $(v_k)_{k \to \infty}$. Therefore, their respective limsup and liminf are equal.

Let's apply this to our bounds. For the right-hand side of $(*)$:
\[
\limsup_{t \to \infty} \frac{N(t)}{a(N(t))} = \limsup_{k \to \infty} \frac{k}{a(k)} = \frac{1}{\liminf_{k \to \infty} \frac{a(k)}{k}} = \frac{1}{\underline{r}}.
\]
For the left-hand side of $(*)$, we rewrite it as $\frac{N(t)}{N(t)+1} \cdot \frac{N(t)+1}{a(N(t)+1)}$. Since $\lim_{t \to \infty} \frac{N(t)}{N(t)+1} = 1$, we have:
\[
\limsup_{t \to \infty} \frac{N(t)}{a(N(t)+1)} = 1 \cdot \limsup_{t \to \infty} \frac{N(t)+1}{a(N(t)+1)} = \limsup_{k \to \infty} \frac{k}{a(k)} = \frac{1}{\underline{r}}.
\]

\textbf{Step 4: Conclusion by the squeeze theorem.}
Taking the $\limsup$ across the inequality $(*)$ yields:
\[
\frac{1}{\underline{r}} \leq \limsup_{t \to \infty} \frac{N(t)}{t} \leq \frac{1}{\underline{r}}.
\]
This proves the first relation: $\limsup_{t \to \infty} \frac{N(t)}{t} = \frac{1}{\underline{r}}$.
The argument for the second relation is identical. Taking the $\liminf$ across $(*)$:
\[
\liminf_{t \to \infty} \frac{N(t)}{a(N(t)+1)} \leq \liminf_{t \to \infty} \frac{N(t)}{t} \leq \liminf_{t \to \infty} \frac{N(t)}{a(N(t))}.
\]
Both the left and right bounds evaluate to $\frac{1}{\limsup_{k \to \infty} \frac{a(k)}{k}} = \frac{1}{\overline{r}}$. Thus:
\[
\frac{1}{\overline{r}} \leq \liminf_{t \to \infty} \frac{N(t)}{t} \leq \frac{1}{\overline{r}}.
\]
This proves the second relation: $\liminf_{t \to \infty} \frac{N(t)}{t} = \frac{1}{\overline{r}}$.
\hfill$\blacksquare$

\textbf{Step 5: System of equations and uniqueness of the limit.}
Injecting Lemma~\ref{lem:density} into the equations from Step 3:
\[
\overline{r} = z + \frac{y - z}{\underline{r}}, \quad \underline{r} = z + \frac{y - z}{\overline{r}}.
\]
Multiplying the first equation by $\underline{r}$, the second by $\overline{r}$:
\[
\overline{r} \underline{r} = z \underline{r} + (y - z), \quad \underline{r} \overline{r} = z \overline{r} + (y - z).
\]
By equality of the right-hand sides:
\[
z \underline{r} + (y - z) = z \overline{r} + (y - z) \implies z \underline{r} = z \overline{r} \implies \underline{r} = \overline{r} \quad (\text{since } z > 0).
\]
The limit $r = \lim_{n \to \infty} \frac{a(n)}{n}$ therefore exists.
\textbf{Step 6: Solving the characteristic equation.}
The equation $r = z + \frac{y - z}{r}$ becomes $r^2 - z r - (y - z) = 0$. Since $r > 0$ (as $a(n) \geq n \min(y,z) + O(1)$), we have:
\[
\boxed{r = \frac{z + \sqrt{z^2 + 4(y - z)}}{2} = r_0.}
\]
\end{proof}

\begin{remark}[Alternative proof via the morphic method]
An alternative proof of Theorem 2.1 can be derived from the morphic nature of the sequences, as established by Fokkink and Joshi \cite{fokkink cloitre xyz.tex}. They show that the characteristic sequence of $S(x,y,z)$ is generated by a morphism whose adjacency matrix (on the alphabet $\{0,1\}$) is $M = \begin{psmallmatrix} z-1 & y-1 \\ 1 & 1 \end{psmallmatrix}$.

The limit $\lim_{n\to\infty} a(n)/n$ corresponds to the Perron-Frobenius eigenvalue of this matrix. The characteristic polynomial is $\lambda^2 - z\lambda - (y-z) = 0$, whose positive root is precisely $r_0$. This method not only confirms the result but also extends its validity to all cases where $y>1$ and $z>1$, including $z>y$.
\end{remark}

\section{Beatty sequences}\label{sec:beatty}

In this section, we show that sequences from the families $S(x,Z+1,Z)$ and $S(x,Z,Z+1)$ correspond to non-homogeneous Beatty sequences. To distinguish these specific integer values from the formal parameters of the general definition, we let $Z \ge 2$ be a fixed integer.
\begin{remark}[Theoretical Context]
The connection between these sequences and Beatty sequences is not a coincidence. The work of Fokkink and Joshi provides a general theoretical framework for these observations. They prove that all hiccup sequences are morphic (see their Theorem 6). More specifically, for the cases where $|y-z|=1$ and the starting value $x$ is not too large ($x \le z$), they show the generating morphism is Sturmian, which has a well-known correspondence with Beatty sequences (see their Theorem 12). The proofs in this section provide direct, constructive demonstrations of this principle for specific, important subfamilies.
\end{remark}

\subsection{The family A: $S(x,Z+1,Z)$}

\begin{myprop}{\texorpdfstring{Closed-form representations for $S(x,Z+1,Z)$}{Closed-form representations for S(x,Z+1,Z)}}{beatty-a}
The irrational slope governing this family is
\[
r_{A}=\frac{Z+\sqrt{Z^{2}+4}}{2}.
\]
For different starting values $x$, the sequences are given by the
following Beatty sequences:
\begin{itemize}
\item For $x=0$:
\[
Z=2\implies a(n)=\Bigl\lfloor(1+\sqrt{2})\,n-\tfrac{3+4\sqrt{2}}{2+\sqrt{2}}\Bigr\rfloor,\;n\ge2.
\]
\[
Z\ge3\implies a(n)=\Bigl\lfloor r_{A}\,n-\tfrac{(Z+1)\,r_{A}-1}{r_{A}+1}\Bigr\rfloor,\;n\ge1.
\]
\item For $x=1,Z\ge2$:
\[
a(n)=\Bigl\lfloor r_{A}\,n-\tfrac{Z\,r_{A}-1}{r_{A}+1}\Bigr\rfloor,\;n\ge1.
\]
\item For $x\ge2,Z\ge2$:
\[
a(n)=\Bigl\lfloor r_{A}\,n-\tfrac{(Z-x)\,r_{A}-1}{r_{A}+1}\Bigr\rfloor,\;n\ge1.
\]
\end{itemize}
For the specific case $x=Z+1$, this simplifies to the quasi-homogeneous
Beatty sequence
\[
a(n)=\lfloor n\,r_{A}+1\rfloor=\lceil n\,r_{A}\rceil.
\]
\end{myprop}

\begin{remark}[Theoretical Justification of Proposition 3.1 via the work of Fokkink et al.]
It is instructive to note that the explicit formulas in Proposition 3.1 can be rigorously justified using the theoretical framework developed by Fokkink and Joshi \cite{fokkink cloitre xyz.tex}. Their approach, while not aimed at deriving general formulas for every case, provides the necessary tools. The proof for a given family, such as $S(x, Z+1, Z)$, proceeds in three steps:

\begin{enumerate}[label=\textbf{Step \arabic*} : , leftmargin=*]
    \item \textbf{Existence of the Beatty form.} Theorem 15 in Fokkink et al. states that a $(0,x,y,z)$-hiccup sequence is a Beatty sequence if $x \le z$ and $|y-z|=1$. For our family, $|(Z+1)-Z|=1$. Thus, for all cases where $x \le Z$, the existence of a Beatty form is guaranteed. Proposition 3.1 is in fact more general, providing formulas that hold for any starting value $x$.

    \item \textbf{Determining the slope.} The slope of the Beatty sequence corresponds to the largest eigenvalue of the adjacency matrix of the generating morphism. For this family, the matrix is $\begin{psmallmatrix} Z-1 & Z \\ 1 & 1 \end{psmallmatrix}$. Its characteristic polynomial is $\lambda^2 - Z\lambda - 1 = 0$, whose positive root is precisely $r_A = \frac{Z+\sqrt{Z^2+4}}{2}$. The slope is thus confirmed.

    \item \textbf{Determining the offset.} Once the slope $r_A$ is known, the offset $\gamma$ in the formula $a(n)=\lfloor n r_A - \gamma \rfloor$ is uniquely determined by the initial condition $a(1)=x$. The condition $\lfloor r_A - \gamma \rfloor = x$ constrains $\gamma$ to an interval of length 1. The formulas in Proposition 3.1 provide the unique value of $\gamma$ that satisfies this constraint for each value of $x$.
\end{enumerate}
This approach shows how the general classification work of Fokkink et al. and the explicit formulas in this article complement and mutually validate each other.
\end{remark}

We present a constructive proof for the special cases where $x=Z+1$, namely for the sequences $S(Z+1, Z+1, Z)$ and $S(Z+1, Z, Z+1)$ (Theorems 3.2 and 3.4). This proof relies directly on the properties of Beatty sequences rather than the more general theory of Sturmian words.

\begin{mythm}{\texorpdfstring{Special case for $x=Z+1$}{Special case for x=Z+1}}{beatty-a-special}
Let $Z\ge2$ be an integer. Define $r_{A}=\frac{Z+\sqrt{Z^{2}+4}}{2}$ and the
sequence
\[
a(n)\;=\;\bigl\lceil n\,r_{A}\bigr\rceil.
\]
Then $a(n)$ satisfies the self-referential definition of the sequence
$S(Z+1,\,Z+1,\,Z)$.
\end{mythm}

\begin{proof}
We must show that $a(n) = \lceil n r_A \rceil$ satisfies the three conditions for $S(Z+1, Z+1, Z)$: (i) initial value, (ii) increment for a "hit", and (iii) increment for a "miss".

\textbf{Step 1: Initial value and possible increments.}
For $Z \ge 2$, we have $Z < r_A < Z+1$, so the initial value is $a(1) = \lceil r_A \rceil = Z+1$, which matches $x$.
The increment $d_k = a(k) - a(k-1) = \lceil k r_A \rceil - \lceil (k-1) r_A \rceil$ is necessarily either $\lfloor r_A \rfloor = Z$ or $\lceil r_A \rceil = Z+1$.

\textbf{Step 2: Linking increments to hits and misses.}
The proof hinges on a standard result from the theory of Beatty sequences: the value of the increment $d_k$ is determined by whether the index $k$ is a member of the sequence's image.
\begin{itemize}[leftmargin=*]
    \item \textbf{Case "hit" ($y=Z+1$):} A "hit" at index $k \ge 2$ means $k \in \{a(1), \dots, a(k-1)\}$. For a Beatty sequence of the form $\lceil n\alpha \rceil$, this condition is known to be equivalent to the increment $d_k = a(k) - a(k-1)$ taking its \emph{larger} possible value. Here, this is $\lceil r_A \rceil = Z+1$. This matches the required increment $y$.

    \item \textbf{Case "miss" ($z=Z$):} A "miss" at index $k$ means $k \notin \{a(1), \dots, a(k-1)\}$. By the Beatty-Rayleigh theorem, such a $k$ must belong to the complementary Beatty sequence. This condition is equivalent to the increment $d_k$ taking its \emph{smaller} possible value. Here, this is $\lfloor r_A \rfloor = Z$. This matches the required increment $z$.
\end{itemize}
Since the initial value and the rules for both types of increments are satisfied, the sequence $a(n)=\lceil n r_A \rceil$ is identical to $S(Z+1, Z+1, Z)$.
\end{proof}

\subsection{The family B: $S(x,Z,Z+1)$}

\begin{myprop}{\texorpdfstring{Closed-form representations for $S(x,Z,Z+1)$}{Closed-form representations for S(x,Z,Z+1)}}{beatty-b}
The irrational slope for this family is
\[
r_{B}=\frac{Z+1+\sqrt{Z^{2}+2Z-3}}{2}.
\]
The sequences are given by:
\begin{itemize}
\item For $x=0,Z\ge2$:
\[
a(n)=\Bigl\lfloor r_{B}\,n-\tfrac{(Z-1)\,r_{B}+1}{r_{B}-1}\Bigr\rfloor,\;n\ge1.
\]
\item For $x=1,Z\ge2$:
\[
a(n)=\Bigl\lfloor r_{B}\,n-\tfrac{(Z-2)\,r_{B}+1}{r_{B}-1}\Bigr\rfloor,\;n\ge2.
\]
\item For $x\ge2,Z\ge2$:
\[
a(n)=\Bigl\lfloor r_{B}\,n-\tfrac{(Z-x)\,r_{B}+1}{r_{B}-1}\Bigr\rfloor,\;n\ge1.
\]
\end{itemize}
For the specific case $x=Z+1$, this simplifies to the quasi-homogeneous
Beatty sequence
\[
a(n)=\lfloor n\,r_{B}+1\rfloor=\lceil n\,r_{B}\rceil.
\]
\end{myprop}

\begin{mythm}{\texorpdfstring{Special case for $x=Z+1$}{Special case for x=Z+1}}{beatty-b-special}
Let $Z\ge2$ be an integer. Define $r_{B}=\frac{Z+1+\sqrt{Z^{2}+2Z-3}}{2}$ and the sequence
\[
a(n)\;=\;\bigl\lceil n\,r_{B}\bigr\rceil.
\]
Then $a(n)$ satisfies the self-referential definition of the sequence $S(Z+1,\,Z,\,Z+1)$.
\end{mythm}

\begin{proof}
The logic is parallel to the previous theorem, but with the roles of the increments reversed.
\textbf{Step 1: Initial value and possible increments.}
For $Z \ge 2$, we have $Z < r_B < Z+1$. The initial value is $a(1) = \lceil r_B \rceil = Z+1$, matching $x$. The increment $d_k = a(k) - a(k-1)$ is necessarily either $\lfloor r_B \rfloor = Z$ or $\lceil r_B \rceil = Z+1$.

\textbf{Step 2: Linking increments to hits and misses.}
\begin{itemize}[leftmargin=*]
    \item \textbf{Case "miss" ($z=Z+1$):} A "miss" at index $k$ means $k \notin \{a(1), \dots, a(k-1)\}$. This condition forces the increment $d_k$ to take its \emph{larger} value, which is $\lceil r_B \rceil = Z+1$. This corresponds to the required increment $z$.

    \item \textbf{Case "hit" ($y=Z$):} A "hit" at index $k$ means $k \in \{a(1), \dots, a(k-1)\}$. This forces the increment $d_k$ to take its \emph{smaller} value, which is $\lfloor r_B \rfloor = Z$. This corresponds to the required increment $y$.
\end{itemize}
The initial value and increment rules match the definition of $S(Z+1, Z, Z+1)$.
\end{proof}

\section{The Case \texorpdfstring{$y=0, z\ge2$}{y=0, z>=2}: Linear Recurrences}\label{sec:recurrence}

When the hit increment is zero, the sequences exhibit ultimately periodic behavior and satisfy linear recurrences.
\begin{mythm}{\texorpdfstring{Linear recurrence for $y=0$}{Linear recurrence for y=0}}{recurrence}
Let $(a(k))$ be $S(x,0,z)$ with $x\in\Z$ and integer $z\ge2$. Then $\lim_{k\to\infty}\frac{a(k)}{k}=z-1$, and for sufficiently large $k$, $a(k)$ satisfies
\[
a(k)-a(k-1)-a(k-z)+a(k-z-1)=0.
\]
\end{mythm}

\begin{proof}
We establish a sequence of properties that build upon each other to prove both claims.
\textbf{Property 1 (Modular invariance):} All terms satisfy $a(j)\equiv x\pmod z$.
\textit{Proof:} By induction. Base case: $a(1)=x\equiv x\pmod z$.
Inductive step: Assume $a(1),\dots,a(k-1)\equiv x\pmod z$.
If $k\in\mathcal{A}_{k-1}$, then $a(k)=a(k-1)+0\equiv x\pmod z$.
If $k\notin\mathcal{A}_{k-1}$, then $a(k)=a(k-1)+z\equiv x\pmod z$.
\textbf{Property 2 (Membership characterization):} For sufficiently large $k$, 
$k\in\mathcal{A}_{k-1}$ if and only if $k\equiv x\pmod z$.
\textit{Proof:} From property 1, $\mathcal{A}_{k-1} \subseteq \{x+jz : j \geq 0\}$. Since $a(n)\to\infty$ and increments are non-negative, for large enough $k$,
all values $x+jz$ with $x+jz < k$ appear in $\mathcal{A}_{k-1}$. Thus for large $k$: $k\in\mathcal{A}_{k-1} \iff k\equiv x\pmod z$.

\textbf{Property 3 (Periodic increments):} For large $k$, the increment sequence becomes periodic with period $z$:
\[
\Delta a(k) = a(k)-a(k-1) = \begin{cases}
0 & \text{if }k\equiv x\pmod z\\
z & \text{otherwise}
\end{cases}
\]
\textit{Proof:} Direct consequence of property 2 and the definition of $S(x,0,z)$.
\textbf{Property 4 (Linear recurrence):} For large $k$, $a(k)$ satisfies
$a(k)-a(k-1)-a(k-z)+a(k-z-1)=0$.
\textit{Proof:} From property 3, $\Delta a(k)$ has period $z$, so $\Delta a(k) = \Delta a(k-z)$. Therefore:
\begin{align*}
a(k)-a(k-1) &= \Delta a(k) = \Delta a(k-z) = a(k-z)-a(k-z-1)
\end{align*}
Rearranging gives the claimed recurrence.

\textbf{Property 5 (Asymptotic slope):} $\lim_{k\to\infty}\frac{a(k)}{k}=z-1$.
\textit{Proof:} In each period of length $z$, exactly one increment is 0 (when $k\equiv x\pmod z$)
and $z-1$ increments are $z$. Thus the average increment per step is:
\[
\frac{1 \cdot 0 + (z-1) \cdot z}{z} = \frac{z(z-1)}{z} = z-1.
\]
Since the increments become exactly periodic, this average becomes the asymptotic slope.
\end{proof}

\section{Vanishing discriminant and combinatorial connections}\label{sec:vanishing}

When the discriminant of the characteristic equation vanishes ($z^2 + 4(y-z) = 0$), the sequences exhibit particularly interesting combinatorial structures. For integer parameters, this occurs when $z=2K$ and $y=2K-K^{2}$. We explore several cases that reveal connections to combinatorics.
\subsection{The sequence $S(3,1,2)$ (OEIS \oeis{A080036})}

This sequence corresponds to the Ramsey core numbers $rc(2,n)$ and is a notable case of a vanishing discriminant ($z^2+4(y-z) = 2^2+4(1-2)=0$).
\begin{mythm}{Ramsey core numbers}{ramsey}
The sequence $S(3,1,2)$ is given by the closed form
\[
a(n)=n+\left\lfloor \frac{\sqrt{8n-7}+3}{2}\right\rfloor.
\]
A key property of this sequence is that its set of values, $\operatorname{im}(a)$, is precisely the set of positive integers that are not of the form $s_j = \frac{j(j+1)}{2} + 1$ for $j \ge 1$ (i.e., triangular numbers plus one).
\end{mythm}

\begin{proof}
Let's define the function $f(n) = n+\left\lfloor \frac{\sqrt{8n-7}+3}{2}\right\rfloor$ and show it satisfies the recurrence for $a(n) = S(3,1,2)$. The proof proceeds by induction, showing that $a(n) = f(n)$ for all $n \ge 1$.
\textbf{Step 1: Base case.}
For $n=1$, the definition gives $a(1)=3$. The formula gives $f(1) = 1 + \lfloor (1+3)/2 \rfloor = 1+2 = 3$. The base case holds.

\textbf{Step 2: Analysis of the increment of $f(n)$.}
Let's analyze the increment $\Delta f(n) = f(n) - f(n-1)$.
\[
\Delta f(n) = 1 + \left( \left\lfloor \frac{\sqrt{8n-7}+3}{2}\right\rfloor - \left\lfloor \frac{\sqrt{8(n-1)-7}+3}{2}\right\rfloor \right).
\]
The difference between the floor terms is either 0 or 1. It is exactly 1 if and only if $n$ is of the form $s_k = \frac{k(k-1)}{2}+1$ for some integer $k \ge 2$. 
Therefore, the increment of $f(n)$ is:
\[
f(n)-f(n-1)=\begin{cases}
2 & \text{if } n = \frac{k(k-1)}{2}+1 \text{ for some } k \ge 2,\\
1 & \text{otherwise}.
\end{cases}
\]

\textbf{Step 3: Inductive equivalence.}
The increment for $a(n)=S(3,1,2)$ is $y=1$ for a "hit" and $z=2$ for a "miss". We must show that the condition for an increment of 2 is the same for both sequences. That is:
\[
n \text{ is a "miss" for } S(3,1,2) \iff n = \frac{k(k-1)}{2}+1.
\]
Assuming the induction hypothesis $a(j)=f(j)$ for $j<n$, a "miss" occurs if $n \notin \{a(1), \dots, a(n-1)\} = \{f(1), \dots, f(n-1)\}$. It is a known property that the set of integers avoided by $f(n)$ is precisely $\{ \frac{k(k-1)}{2}+1 \mid k \ge 2 \}$.
Thus, the increments match at every step, and so $a(n)=f(n)$ for all $n$.
\end{proof}

\begin{figure}[H]
\centering
\begin{tikzpicture}[scale=0.8]
    \draw[gray, thin] (-1,-1) grid (4,4);
    \foreach \y in {0,...,3} {
        \foreach \x in {0,...,\y} {
            \pgfmathtruncatemacro{\idx}{\y*(\y+1)/2 + \x + 1}
            \node at (\x,3-\y) [circle, fill=blue!60, inner sep=2pt] {};
            \node at (\x,3-\y) [above right] {\scriptsize \idx};
        }
    }
    \foreach \y in {0,...,3} {
        \pgfmathtruncatemacro{\idx}{\y*(\y+1)/2 + 1}
        \node at (0,3-\y) [circle, fill=red!80, inner sep=2pt] {};
    }
    \node[align=left] at (7,2.5) {\small\color{red!80}Red: avoided points\\ \small\color{blue!60}Blue: sequence points};
    \node at (7,1) {\small $s_2=2, s_3=4, s_4=7, s_5=11$};
\end{tikzpicture}
\caption{Triangular grid interpretation of $S(3,1,2)$}
\label{fig:triangular}
\end{figure}
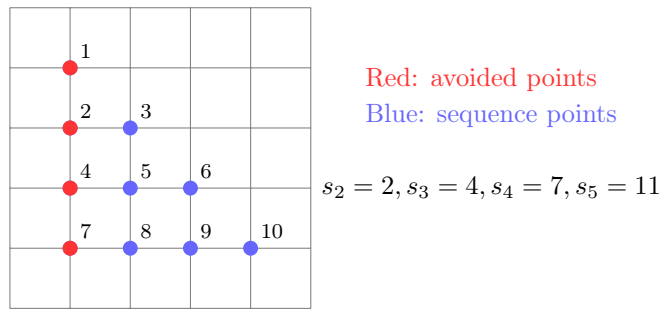

\subsection{The sequence $S(4,1,2)$ (OEIS \oeis{A217334})}

This sequence, also corresponding to a vanishing discriminant case, has a powerful combinatorial interpretation related to covering a grid with squares.
\begin{mythm}{Thumbtack sequence}{thumbtacks}
The sequence $S(4,1,2)$ has the closed form
\[
a(n)=n+\floor{\sqrt{4n-3}}+2.
\]
Combinatorially, it represents the minimum number of vertices required to define $n$ unit squares placed in an outward spiral on the $\mathbb{Z}^2$ grid, starting from the origin.
\end{mythm}

\begin{proof}
Let $g(n) = n+\floor{\sqrt{4n-3}}+2$. We prove by induction that $g(n)$ is identical to $a(n)$ generated by the $S(4,1,2)$ rule.
\textbf{Step 1: Base case.}
For $n=1$, the rule gives $a(1)=4$. The formula gives $g(1) = 1+\floor{\sqrt{1}}+2 = 4$. The base case holds.

\textbf{Step 2: Analysis of the increment of $g(n)$.}
The increment $\Delta g(n) = g(n)-g(n-1)$ is:
\[
\Delta g(n) = 1 + \left( \floor{\sqrt{4n-3}} - \floor{\sqrt{4(n-1)-3}} \right).
\]
The difference of the floor terms is 1 (making the total increment 2) if and only if $n$ is a value of the form $\lfloor k^2/4 \rfloor + 1$. These points correspond to the "misses". For all other values of $n$, the increment is 1.

\textbf{Step 3: Inductive equivalence.}
The increment of $a(n)$ is 2 if $n$ is a "miss" and 1 if $n$ is a "hit". Assuming $a(k)=g(k)$ for $k<n$, we must show that $n$ is a miss for $a(n)$ if and only if $n$ is a "jump point" for $g(n)$. This equivalence can be established by showing that the set of integers not in the image of $g$ corresponds to the set of jump points. The geometric interpretation confirms this: a "miss" (increment 2) corresponds to starting a new corner of the spiral, which happens at these jump points. A "hit" (increment 1) corresponds to adding a square along an edge. Since the initial values and all subsequent increments match, $a(n)=g(n)$ for all $n \ge 1$.
\end{proof}

\begin{figure}[H]
\centering
\begin{subfigure}{0.48\textwidth}
    \centering
    \begin{tikzpicture}[scale=0.5]
        \draw[gray!30, thin] (-2.5,-2.5) grid (2.5,2.5);
        \fill[blue!20] (0,0) rectangle (1,1);
        \foreach \p in {(0,0), (1,0), (0,1), (1,1)}
            \node[red, cross out, thick, draw, inner sep=1pt] at \p {};
    \end{tikzpicture}
    \caption{$n=1$: 4 new vertices. $a(1)=4$.}
\end{subfigure}
\hfill
\begin{subfigure}{0.48\textwidth}
    \centering
    \begin{tikzpicture}[scale=0.5]
        \draw[gray!30, thin] (-2.5,-2.5) grid (2.5,2.5);
        \fill[gray!10] (0,0) rectangle (1,1);
        \fill[blue!20] (1,0) rectangle (2,1);
        \foreach \p in {(0,0), (1,0), (0,1), (1,1)}
            \node[gray!80, cross out, thick, draw, inner sep=1pt] at \p {};
        \node[red, cross out, thick, draw, inner sep=1pt] at (2,0) {};
        \node[red, cross out, thick, draw, inner sep=1pt] at (2,1) {};
    \end{tikzpicture}
    \caption{$n=2$: 2 new vertices. $a(2)=6$.}
\end{subfigure}
\vspace{0.8cm}
\begin{subfigure}{0.48\textwidth}
    \centering
    \begin{tikzpicture}[scale=0.5]
        \draw[gray!30, thin] (-2.5,-2.5) grid (2.5,2.5);
        \fill[gray!10] (0,0) rectangle (1,1);
        \fill[gray!10] (1,0) rectangle (2,1);
        \fill[blue!20] (1,1) rectangle (2,2);
        \foreach \p in {(0,0), (1,0), (2,0), (0,1), (1,1), (2,1)}
            \node[gray!80, cross out, thick, draw, inner sep=1pt] at \p {};
        \node[red, cross out, thick, draw, inner sep=1pt] at (1,2) {};
        \node[red, cross out, thick, draw, inner sep=1pt] at (2,2) {};
    \end{tikzpicture}
    \caption{$n=3$: 2 new vertices. $a(3)=8$.}
\end{subfigure}
\hfill
\begin{subfigure}{0.48\textwidth}
    \centering
    \begin{tikzpicture}[scale=0.5]
        \draw[gray!30, thin] (-2.5,-2.5) grid (2.5,2.5);
        \fill[gray!10] (0,0) rectangle (1,1);
        \fill[gray!10] (1,0) rectangle (2,1);
        \fill[gray!10] (1,1) rectangle (2,2);
        \fill[blue!20] (0,1) rectangle (1,2);
        \foreach \p in {(0,0), (1,0), (2,0), (0,1), (1,1), (2,1), (1,2), (2,2)}
            \node[gray!80, cross out, thick, draw, inner sep=1pt] at \p {};
        \node[red, cross out, thick, draw, inner sep=1pt] at (0,2) {};
    \end{tikzpicture}
    \caption{$n=4$: 1 new vertex. $a(4)=9$.}
\end{subfigure}
\caption{Spiral covering interpretation of $S(4,1,2)$.}
\label{fig:spiral}
\end{figure}
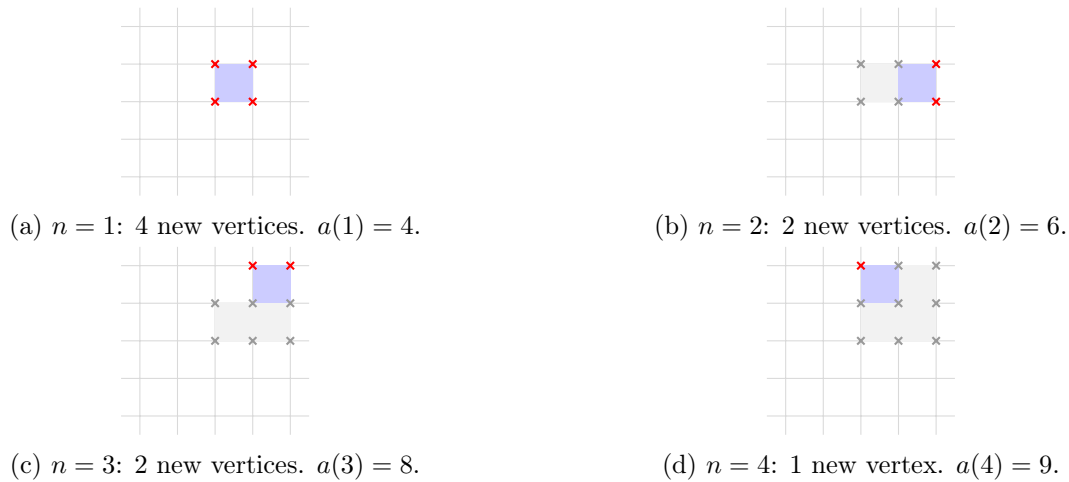

\subsection{The sequence $S(5,1,2)$ (OEIS \oeis{A080353})}
\label{subsec:S512}

This case reveals a surprisingly deep and elegant connection to hexagonal lattices. While its closed form is complex, it precisely captures the recursive "hit/miss" dynamic.
\begin{mythm}{Hexagonal lattice sequence}{hexagonal}
The sequence $S(5,1,2)$ is given by the closed form:
\[
a(n) = n + 3m + 1 + \left\lfloor \frac{n - 1 - \frac{3m(m-1)}{2}}{m} \right\rfloor,
\]
where the parameter $m$, corresponding to the layer of a hexagonal spiral, is given by:
\[
m = \left\lceil \frac{-1 + \sqrt{1 + \frac{8n}{3}}}{2} \right\rceil.
\]
\end{mythm}

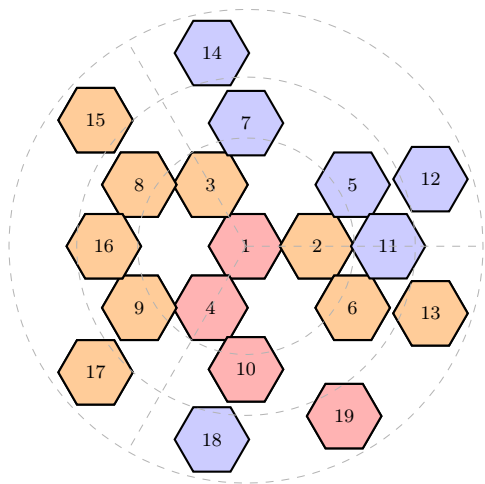
\begin{figure}[H]
\centering
\begin{tikzpicture}[
    scale=0.9, every node/.style={transform shape}
]
    \tikzset{
      hex/.style       = {regular polygon, regular polygon sides=6, minimum size=1.1cm, draw, thick, inner sep=1pt, font=\scriptsize},
      guideline/.style = {dashed, gray!60, very thin},
      primary/.style   = {hex, fill=red!30},
      secondary/.style = {hex, fill=orange!40},
      hit/.style       = {hex, fill=blue!20}
    }
    \coordinate (1) at (0,0);
    \coordinate (2) at (0:1.05); \coordinate (3) at (120:1.05); \coordinate (4) at (240:1.05);
    \coordinate (5) at (30:1.82); \coordinate (6) at (-30:1.82);
    \coordinate (7) at (90:1.82);
    \coordinate (8) at (150:1.82); \coordinate (9) at (210:1.82); \coordinate (10) at (270:1.82);
    \coordinate (11) at (0:2.1);
    \coordinate (12) at (20:2.9); \coordinate (13) at (-20:2.9);
    \coordinate (14) at (100:2.9); \coordinate (15) at (140:2.9); \coordinate (16) at (180:2.1);
    \coordinate (17) at (220:2.9); \coordinate (18) at (260:2.9); \coordinate (19) at (300:2.9);
    \def\primarymisses{1,4,10,19}
    \def\secondarymisses{2,3,6,8,9,13,15,16,17}
    \foreach \n in {1,...,19}{
        \node[hit] at (\n) {\n};
        \foreach \s in \secondarymisses {
            \ifnum\n=\s \node[secondary] at (\n) {\n}; \fi
        }
        \foreach \p in \primarymisses {
            \ifnum\n=\p \node[primary] at (\n) {\n}; \fi
        }
    }
    \draw[guideline] (0,0) -- (0:3.5);
    \draw[guideline] (0,0) -- (120:3.5);
    \draw[guideline] (0,0) -- (240:3.5);
    \draw[guideline, rotate=60] (0,0) circle (1.6);
    \draw[guideline, rotate=60] (0,0) circle (2.5);
    \draw[guideline, rotate=60] (0,0) circle (3.5);
\end{tikzpicture}
\caption{Structure of the hexagonal "snowflake" generated by $S(5,1,2)$ (up to $n=19$). Misses (red and orange) appear on outer edges.}
\label{fig:hex}
\end{figure}

\begin{proof}[Proof outline and geometric interpretation]
A full algebraic proof showing that this closed form satisfies the $S(5,1,2)$ recurrence is technical and omitted for brevity. The logic relies on demonstrating that the increment of the closed form, $\Delta a(n) = a(n)-a(n-1)$, is equal to 2 if and only if the index $n$ is a "miss" for the recursive sequence. This set of "misses" has a rich geometric structure tied to the hexagonal lattice.
The misses can be decomposed into two classes, as visualized in Figure~\ref{fig:hex}:
\begin{itemize}
    \item \textbf{Primary misses}: Indices of the form $H_k = \frac{3k(k-1)}{2} + 1$ for $k \geq 1$ (e.g., $1, 4, 10, 19, \dots$). These correspond to vertices that initiate new hexagonal layers or shells in the spiral growth.
    \item \textbf{Secondary misses}: Additional indices (e.g., $2, 3, 6, 8, \dots$). These correspond to other "corner" or "edge" vertices required to complete a layer before the next primary miss is reached.
\end{itemize}
The closed form perfectly models this complex distribution of miss events, confirming the deep connection between the simple $S(5,1,2)$ recurrence and the geometry of hexagonal lattices.
\end{proof}

\section{Connections to deeper combinatorics and algorithm analysis}\label{sec:metafib}

Beyond the connections to classical number theory (Beatty sequences) and geometric lattices, the $S(x,y,z)$ family reveals deep ties to other areas of modern combinatorics, notably the study of tree structures, the combinatorics of words, and the analysis of algorithms.
\subsection{Connection to trees and Meta-Fibonacci sequences}

The connection between our $S(x,y,z)$ sequences and meta-Fibonacci sequences can be made precise through the combinatorial work of Ruskey and Deugau~\cite{Ruskey}. Their results establish a formal bridge between a subfamily of our sequences and a well-studied class of meta-Fibonacci recurrences, revealing a shared combinatorial foundation related to k-ary trees.
\begin{mythmruskey}{}{ruskey-deugau}
For all $n \geq 1$ and $k \geq 2$, the sequence $b_k(n) = S(k+1,1,k+1)$ is given by the explicit, non-recursive formula:
\[
b_k(n) = n + k \cdot a_{0,k}(n),
\]
where $a_{0,k}(n)$ counts the number of leaves in a forest of complete k-ary trees with a total of $n$ nodes.
\end{mythmruskey}

This identity is powerful, as it reveals that the "hit/miss" rule for this subfamily is governed by the number of leaves in growing k-ary trees. For example, for $k=2$, the sequence $S(3,1,3)$ is given by \oeis{A045412}. This sequence, in turn, has a further, surprising connection to computer science. As noted by C. J. Quines, it precisely describes the indices of all terms greater than 1 in the sequence \oeis{A182105}, which counts the number of elements being merged at each step of a bottom-up merge sort algorithm. Remarkably, the relevance of this sequence extends to modern AI research, as it appears in recent work on strategy scheduling for automated theorem provers~\cite{Bartek2024}.
\subsection{Connection to Combinatorics on Words}

Other choices of parameters yield sequences that are best understood through the lens of combinatorics on words, where the sequence gives the positions of certain letters in an infinite word generated by a morphism. A prime example is the sequence $S(1,2,4)$, which is identical to \oeis{A080580}. As proven on its OEIS page, this sequence gives the positions of the letter '0' in the infinite word generated as the fixed point of the morphism:
\[
\sigma: \begin{cases}
0 \to 01 \\
1 \to 1101
\end{cases}
\]
Starting with '0' and repeatedly applying $\sigma$ gives:
\[
0 \to 01 \to 011101 \to 0111011101110101 \to \dots
\]
The positions of the '0's in this infinite word are $1, 5, 9, 13, 15, \dots$, which exactly match the sequence $S(1,2,4)$. This shows that the simple hit/miss recurrence can act as a compact and efficient generator for otherwise complex morphic sequences.
\begin{remark}
Fokkink and Joshi \cite{fokkink cloitre xyz.tex} have shown this is a general property: the characteristic sequence of any $S(x,y,z)$ is morphic.
\end{remark}

\section{Concluding Remarks}

The $S(x,y,z)$ sequences provide a rich framework connecting recursive definitions, asymptotic analysis, and combinatorial structures.
Key results include:

\begin{itemize}[leftmargin=*]
\item \textbf{Linear asymptotics} for $\Delta>0$ and $y>z$ with slope $r_{0}$.
\item \textbf{Beatty sequence generation} for $S(x,Z+1,Z)$ and $S(x,Z,Z+1)$ when $Z\geq2$.
\item \textbf{Periodic increments} and linear recurrences for $y=0$, $z\ge2$.
\item \textbf{Combinatorial interpretations for discriminant-zero cases:} When the characteristic discriminant vanishes, the sequences generate rich combinatorial structures tied to the geometry of regular lattices.
The emergence of patterns based on fundamental \textbf{triangular} ($S(3,1,2)$), \textbf{square} ($S(4,1,2)$), and \textbf{hexagonal} ($S(5,1,2)$) tessellations suggests potential ramifications in mathematical crystallography.
\end{itemize}

\begin{tcolorbox}[colback=orange!10,colframe=orange!80,title={\textbf{Open Problems}}]
\begin{enumerate}
\item Prove that $a(n) = r_0 n + O(1)$ in certain cases where $y > z$.
\item Establish a combinatorial proof that Ramsey core numbers $rc(2,n) = S(3,1,2)$ by proving the equivalence between the combinatorial definition of $rc(2,n)$ and our triangular grid interpretation.
\item Study the class of sequences $S(x,y(n),z(n))$ where $y$ and $z$ are bounded or unbounded integer sequences.
\item Investigate the 4-parameter class $S(x,y,z,t)$ where increments depend on membership in three disjoint sets partitioning $\N$.
\end{enumerate}
\end{tcolorbox}

\appendix
\section{OEIS Compendium of Hiccup Sequences}
The following table, adapted and extended from the work of Fokkink and Joshi, provides a comprehensive summary of hiccup sequences.
It includes both the sequences officially recognized as such in the OEIS at the time of their paper's writing, and others identified by the authors as belonging to the same family.
The table uses the `(j,x,y,z)`-hiccup sequence notation from Fokkink and Joshi. The family `S(x,y,z)` studied in this article corresponds precisely to the case where `j=0`. In the general case, the parameter `j` modifies the self-referential condition to become $k-j \in \mathcal{A}_{k-1}$.

\begin{table}[H] 
\centering
\caption{Hiccup sequences and their parameters $(j,x,y,z)$.}
\label{tab:oeis_compendium}
\begin{tabular}{lcccc}
\toprule
\textbf{OEIS Number} & \textbf{j} & \textbf{x} & \textbf{y} & \textbf{z} \\
\midrule
\oeis{A000201} & 1 & 1 & 2 & 1 \\
\oeis{A003156} & 1 & 1 & 3 & 1 \\
\oeis{A004956} & 0 & 2 & 2 & 1 \\
\oeis{A007066} & 0 & 1 & 2 & 3 \\
\oeis{A026352} & 1 & 1 & 2 & 3 \\
\oeis{A026356} & 0 & 2 & 2 & 3 \\
\oeis{A045412} & 0 & 3 & 1 & 3 \\
\oeis{A064437} & 0 & 1 & 3 & 2 \\
\oeis{A080578} & 0 & 1 & 1 & 3 \\
\oeis{A080579} & 0 & 1 & 1 & 4 \\
\oeis{A080580} & 0 & 1 & 2 & 4 \\
\oeis{A080590} & 0 & 1 & 3 & 4 \\
\oeis{A080600} & 0 & 4 & 4 & 3 \\
\oeis{A080652} & 0 & 2 & 3 & 2 \\
\oeis{A080667} & 0 & 3 & 4 & 3 \\
\oeis{A080903} & 0 & 1 & 4 & 2 \\
\oeis{A081834} & 0 & 1 & 4 & 3 \\
\oeis{A081835} & 0 & 1 & 5 & 4 \\
\oeis{A081839} & 0 & 0 & 4 & 5 \\
\oeis{A081840} & 0 & 0 & 3 & 4 \\
\oeis{A081841} & 0 & 0 & 3 & 2 \\
\oeis{A081842} & 0 & 0 & 4 & 3 \\
\oeis{A081843} & 0 & 0 & 5 & 4 \\
\oeis{A086377} & 1 & 1 & 3 & 2 \\
\oeis{A086398} & 1 & 1 & 4 & 2 \\
\oeis{A284753} & 0 & 2 & 4 & 2 \\
\bottomrule
\end{tabular}
\end{table}

\end{document}